\documentclass[12pt,a4paper]{amsart}
\usepackage{graphicx,amssymb}
\input xy
\xyoption{all}
\usepackage[all]{xy}
\usepackage{hyperref}

\newtheorem{thm}{Theorem}[section]
\newtheorem{cor}[thm]{Corollary}
\newtheorem{lem}[thm]{Lemma}
\newtheorem{prop}[thm]{Proposition}
\theoremstyle{definition}

\newtheorem{rem}[thm]{Remark}

\numberwithin{equation}{section}
\newcommand{\ZZ}{\mathbb Z}

\newcommand{\PP}{\mathbb P}

\newcommand{\lra}{\longrightarrow}
\newcommand{\hra}{\hookrightarrow}
\newcommand{\ra}{\rightarrow}

\newcommand{\cA}{\mathcal{A}}

\newcommand{\cT}{\mathcal{T}}
\newcommand{\tC}{\tilde{C}}

\newcommand{\cO}{\mathcal{O}}
\newcommand{\cR}{\mathcal{R}}

\DeclareMathOperator{\Aut}{{Aut}}
\DeclareMathOperator{\Hom}{{Hom}}
 
 \DeclareMathOperator{\Ker}{Ker}
\DeclareMathOperator{\Pic}{Pic}

 \DeclareMathOperator{\Ima}{{Im}}
 \DeclareMathOperator{\Min}{Min}
 \DeclareMathOperator{\Cliff}{Cliff}

 \DeclareMathOperator{\im}{Im}

\begin{document}

\title[ ]{Prym varieties of cyclic coverings}%
\author{H. Lange and A. Ortega}

\address{H. Lange\\Mathematisches Institut, Universit\"at Erlangen-N\"urnberg\\Germany}
\email{lange@mi.uni-erlangen.de}
              
\address{A. Ortega \\ Fachbereich Mathematik, Universit\"at Duisburg-Essen\\ Germany}
\email{angela.ortega@uni-due.de}

\thanks{We would like to thank Edoardo Sernesi for some valuable hints concerning the proof of Proposition 4.1. }

%\date{ }%
%\dedicatory{ }%
%\commby{ }%
% ----------------------------------------------------------------
\begin{abstract} The Prym map of type $(g,n,r)$ associates to every cyclic covering of degree $n$ of a curve of 
genus $g$, ramified at a reduced divisor of degree $r$, the corresponding Prym variety. We show that the 
corresponding map of moduli spaces is generically finite in most cases. From this we deduce the dimension of the image of the Prym map.

\end{abstract}

\maketitle

\section{Introduction}

Let $C$ be a smooth projective curve of genus $g \geq 2$ over an algebraically closed field $k$. 
For simplicity we assume the characteristic of $k$ to be 0, although some of the results are valid in greater generality.
To every \'etale double covering $f: \tC \ra C$ one can associate a principally polarized abelian variety 
of dimension $g-1$,
the Prym variety of $f$. The corresponding morphism  
$$
Pr: \cR_g \lra \cA_{g-1}
$$
from the moduli space of \'etale double coverings of curves of genus $g$ to the moduli space of 
principally polarized abelian varieties of dimension $g-1$ is called the Prym map. 
Its properties have been extensively studied. For example it is well known that for $g\geq 6$ it is 
generically injective (see \cite{fs}, \cite{k}) and not injective (see \cite{d}). Moreover 
its infinitesimal properties have been studied as well (see \cite{Be}, \cite{ls}).

In this paper we look at the analogue of the Prym map in the case of cyclic coverings of arbitrary degree. To any 
cyclic covering $f: \tC \ra C$ of degree $n$, branched along a reduced divisor $B$ on C of degree $r$, one can associate 
an abelian variety $P = P(f)$ of dimension $p$, called the {\it Prym variety of the covering} $f$ and defined as the connected component containg 0
of the norm map $J\tC \ra JC$. The canonical principal polarization of $J\tC$ induces a 
polarization on $P$ whose type $D$ depends only on the topological structure of the covering $f$.
  
Giving a such covering is equivalent to give the data $(C, B, \eta)$ where $\eta $ is a 
line bundle over $C$ satisfying $\eta^{\otimes n} \simeq \cO_C(B)$. 
We denote by $\cR_g(n,r) $ the corresponding coarse moduli space. If $\cA _{p,D}$ denotes 
the moduli space of polarised abelian varieties of dimension $p$ and polarisation $D$, the morphism 
$$
Pr_g(n,r): \cR_g(n,r)  \lra \cA _{p,D}
$$
which associates to a $n$-cyclic  covering its polarized Prym variety is called the {\it Prym map of type} $(g,n,r)$.

In this paper we investigate the differential of the Prym map at general $n$-cyclic branched 
coverings. We will show that in most cases it is injective 
(see Propositions \ref{prop5.1}, \ref{prop5.2}, \ref{prop5.6} and \ref{prop5.7}). 
As a consequence we get that the Prym map itself is generically finite in most cases (see Proposition \ref{prop5.5} 
and Corollary \ref{cor5.8}). From this we can compute the dimension of the image of the Prym map. 
 
In Section 2 we collect some results concerning the ampleness of the twisted canonical line bundle 
$\omega_C \otimes \eta$ which will be used throughout the article. 
Section  3 is dedicated to the study of the injectivity of the Abel-Prym map.  
In Section 4 we prove that the codifferential of the Prym map can be identified with a certain multiplication map of sections, 
which give us a necessary condition to  the codifferential map be surjective at a point $(C,B, \eta)$  of $\cR_g(n,r)$. 
Finally, the consequences of this result are worked out in Section  5.

\section{Twisted canonical line bundles}

Let $C$ be a smooth projective curve of genus $g \geq 2$ over an algebraically closed field $k$ 
(of characteristic 0) and $\eta$ be a line bundle of degree 
$s \geq 0$. In this section we collect some well-known results concerning the line bundle 
$$
L = \omega_C \otimes \eta,
$$ which 
will be applied later.

\begin{lem} \label{lem2.1}
{\em (1)} $L$ is globally generated if and only if either
\begin{itemize}
\item $s \geq 2$ or
\item$s = 1$ and $\eta \neq \cO_C(x)$ for some point $x$ of $C$ or
\item $s = 0$ and $\eta \neq \cO_C(x - y)$ for points $x \neq y$ of $C$.
\end{itemize}
{\em (2)} $L$ is very ample if and only if either
\begin{itemize}
\item $s \geq 3$ or 
\item $s = 2$ and $\eta \neq \cO_C(x + y)$ for points $x$ and $y$ of $C$ or
\item $s=1$ and $\eta \neq \cO_C(x+y - u)$ for points $x, y$ and $u$ of $C$ or
\item $s=0$ and $\eta \neq \cO_C(x+y - u - v)$ for points $x, y, u$ and $v$ of $C$ with $x+y \neq u+v$.
\end{itemize}
\end{lem}
\begin{proof} This is a consequence of Riemann-Roch.
\end{proof}
We need the following immediate corollary of Lemma \ref{lem2.1} (2) in Section 5.
\begin{cor} \label{cor2.2}
If $\deg \eta = 0$, then $\omega_C \otimes \eta$ is very ample
if and only if $\omega_C \otimes \eta^{-1}$ is very ample.
\end{cor}
 Recall that the Clifford index of the curve $C$ is defined by
 $$
 \Cliff (C) = \Min \{\deg M - 2(h^0(M) -1)\;|\; M \in \Pic(C), h^i(M) \geq 2\;  \mbox{for} \; i = 0,1 \}.
 $$
\begin{cor} \label{cor2.3}
Suppose 
$$
\eta^n = \cO_C(B)
$$
with $n \geq 2$ and a reduced effective divisor $B$ of degree $r = ns$ on $C$. Then\\
{\em (1)} Suppose $s = 0$ or $1$. If $g \geq n+1$ and $\Cliff(C) \geq n-1$, then $L$ is globally generated.\\
{\em (2)} Suppose $s = 0, 1$ or $2$. If $g \geq 2n+1$ and $\Cliff(C) \geq 2n-1$, then $L$ is very ample.
\end{cor}

The Clifford index of a general curve equals $[\frac{g-1}{2}]$. Hence for a general curve $C$ of genus $g \geq 2n-1$,
respectively $\geq 4n-1$, every line bundle $L$ of this type is globally generated, respectively very ample. 
\begin{proof}
(1): Let $s = 0$. If $L$ is not globally generated, then $ \eta = \cO(x-y)$ for points $x \neq y$ of $C$ 
according to Lemma \ref{lem2.1} (1). So $nx  \sim ny$ implying that $C$ admits a $g_n^1$. 
The assumption on the genus implies that
it contributes to the Clifford index, which means $\Cliff(C) \leq n-2$.
The proof in the case $s =1$ is the same.\\
(2): Let $s = 0$. If $L$ is not very ample, then $\eta = \cO(x+y-u-v)$ for points $x,y,u,v$ of $C$ with $x+y \neq u+v$
according to Lemma \ref{lem2.1} (2).
So $nx + ny \sim nu + nv$ implying that $C$ admits a $g_{2n}^1$. Again the assumption on the genus implies that it
contributes to the Clifford index, which means that $\Cliff(C) \leq 2n-2$. The proof of the other cases is the same. 
\end{proof}

\begin{cor} \label{cor2.4}
Let $\eta \in \Pic^0(C)$ with $\eta^n = \cO_C$ as above. If $\Cliff(C) \geq 2n-1$, then
the canonical map
$$
H^0(\omega_C \otimes \eta) \otimes H^0(\omega_C \otimes \eta^{-1}) \lra H^0(\omega_C^2)
$$
is surjective. 
\end{cor}
\begin{proof}
According to a classical Theorem of Meis, $\Cliff(C) \leq \frac{g-1}{2}$. So the assumption on the Clifford index implies
$ g \geq 4n-1$ and $\omega_C \otimes \eta$ and $\omega_C \otimes \eta^{-1}$ are very ample by Corollary \ref{cor2.3}.
Hence the asumptions of \cite[Theorem 1]{b} are satisfied, which gives the assertion.
\end{proof}

\section{The Abel-Prym map}

Let $C$ be a smooth projective curve of genus $g \geq 1$ and  
$$
f: \tilde{C} \ra C
$$ 
denote a cyclic covering 
of degree $n$ branched over a reduced divisor $B = x_1 + \cdots + x_r$ on $C$ of degree $r \geq 0$. 
So $\tilde{C}$ is given by a line bundle $\eta$ of degree $s = \frac{r}{n}$ on $C$ satisfying
$$
\eta^n = \cO_C(B)
$$
and we have
$$
f_*\cO_{\tilde{C}} = \bigoplus_{i=0}^{n-1} \eta^{-i}, \qquad  \qquad \omega_{\tilde{C}} = f^*(\omega_C \otimes \eta^{n-1}).
$$
For such a covering $f$ let 
$$
\sigma : \tC \rightarrow \tC
$$ 
denote a map generating the cyclic group of covering maps. 
The automorphism $\sigma$ induces an automorphism on the Jacobian $J\tC$ which we 
denote by $\sigma$ as well.

The Prym variety $P = P(f)$ of the covering $f$ is defined as 
$$
P = \Ima(1-\sigma) = \Ker (1 + \sigma + \cdots + \sigma^{n-1})^0.
$$
Fix a base point $c \in \tilde{C}$ and let $\alpha_c: \tilde{C} \ra J(\tilde{C})$ be the Abel map with respect to $c$,
i.e. $\alpha_c(p) = \cO_{\tilde{C}}(p - c)$. The {\it Abel-Prym map} of the covering $f$ 
is defined as the composition
$$
\pi = \pi_c: \tilde{C} \stackrel{\alpha_c}{\lra} J(\tilde{C}) \stackrel{1 - \sigma}{\lra} P.
$$

\begin{prop}  \label{prop3.1}
{\em (1)} Suppose $\tC$ is not hyperelliptic or $\tC$ is hyperelliptic and $n \geq 3$. Then, for distinct points $p,q \in \tC$,  $\pi(p)=\pi(q)$  if and only if $f$ is 
ramified in $p$ and $q$ and all ramification points have the same image. In particular,  $\pi:  \tC \rightarrow
P$ is injective if $f$ is \'etale.\\
{\em (2)} If $\tC$ is hyperelliptic and $n=2$, then $\pi: \tC \rightarrow P$ is of degree 2 onto its image. 
All ramification points of $f$ have the same image under $\pi$.
\end{prop}
\begin{proof}
By definition of the Abel-Prym map, $\pi(p)=\pi(q)$  if and only if 
$ (1-\sigma)(p-c) \sim (1-\sigma)(q-c)$. This is the case if and only if 
\begin{equation} \label{eq3.1}
p + \sigma (q) \sim q+ \sigma(p).
\end{equation}
If $\tC$ is not hyperelliptic this means  that $p + \sigma (q) = q +\sigma(p)$.  Since $p\neq q$,  this is the case 
if and only if $p$ and $q$ are
ramification points of $f$. \\
Let $\tC$ be a hyperelliptic curve, $\iota$ the hyperelliptic involution, and $p,q \in \tC$ 
distinct points satisfying \eqref{eq3.1}. Again all ramification points have the same image.
If $p\neq q $ are not ramification points, \eqref{eq3.1} implies $\sigma(q) = \iota (p)$ and $\sigma(p)= \iota (q)$. Then
$$
\sigma^2(q) = \sigma \iota(p) = \iota \sigma (p) =\iota^2(q) = q,
$$
since $\sigma$ commutes with the hyperelliptic involution. Similarly, we get $\sigma^2(p)=p$. Then $f$ ramifies
in $p$ and $q$ for $n \geq 3$, a contradiction. The proof in the case $n=2$ is the same as for \cite[Proposition 12.5.2]{bl}
where it was proved for $s \leq 1$. The proof of the fact that also in the hyperelliptic case all ramification points 
of $f$ have the same image under $\pi$ is the same as above.
\end{proof}

In order to describe the differential of the Abel-Prym map, let
$\chi$ denote the generating character of the Galois group $G = \; <\sigma>$ of $f$, i.e. $\chi(\sigma) = \zeta_n$, 
a primitive $n$-th root of unity.
The decomposition of $H^0(\omega_{\tilde{C}})$ into eigenspaces is
\begin{equation} \label{dec}
H^0(\tC, \omega_{\tC}) = \bigoplus_{i=0}^{n-1} H^0(C, \omega_C \otimes \eta^{n-i}),
\end{equation}
where $H^0(\omega_C \otimes \eta^{n-i})$ is the eigenspace of $\chi^i$.
Recall that 
$$
J(\tilde{C}) = H^0(\tC, \omega_{\tilde{C}})^*/H_1(\tilde{C},\ZZ).
$$
In these terms the induced action of $\sigma$ on $H^0(\tC, \omega_{\tilde{C}})^*$ and 
$H_1(\tilde{C},\ZZ)$ is just the analytic,  
respectively rational, representation of $\sigma$. 
We denote 
\begin{equation} \label{eq3.2}
H^0(\tC, \omega_{\tilde{C}})^+ = H^0(C, \omega_C(B)) \quad \mbox{and} \quad 
H^0(\tC, \omega_{\tilde{C}})^- = \bigoplus_{i=1}^{n-1} H^0(C, \omega_C \otimes \eta^{n-i}).
\end{equation} 
Notice that $H^0(\tC, \omega_{\tilde{C}})^+$ is the eigenspace with eigenvalue $1$ of the action of $G$ 
on $H^0(\tC, \omega_{\tilde{C}})$ and $H^0(\tC, \omega_{\tilde{C}})^-$ its complement, 
that is the sum of the other eigenspaces. 
If we define similarly $H_1(\tilde{C},\ZZ)^+$ and $H_1(\tilde{C},\ZZ)^-$, we have
\begin{equation} \label{eq3.3}
P = (H^0(\tC, \omega_{\tilde{C}})^-)^*/H_1(\tilde{C},\ZZ)^-.
\end{equation}
Hence the tangent bundle of $P$ is the trivial bundle 
$$
\cT_P = P \times (H^0(\tC, \omega_{\tilde{C}})^-)^*.
$$
The projectivized differential of the Abel-Prym map $\pi_c$ is by definition the projectivization of the composed map
$$
\cT_{\tC} \stackrel{d\pi_c}{\lra} \cT_P = P \times (H^0(\tC, \omega_{\tilde{C}})^-)^* \lra (H^0(\tC, \omega_{\tilde{C}})^-)^*.
$$
It is a priori a rational map $\tC=P(\cT_{\tC}) \dasharrow P((H^0(\tC, \omega_{\tilde{C}})^-)^*)$.
Let $q_i: (H^0(\tC, \omega_{\tilde{C}})^-)^* \ra H^0(C, \omega_C \otimes \eta^i)^*$ denote the natural 
projection and $P(q_i)$ its projectivization. 
 
\begin{lem} \label{lem3.2}
The composition of the projectivized differential of the Abel-Prym map $\pi_c$ with the map 
$P(q_i)$ for $ 1 \leq i \leq n-1$ is the composed map
$$
\varphi_{\omega_C \otimes \eta^i} \circ f: \tilde{C} \ra C \ra P(H^0(C, \omega_C \otimes \eta^i)^*)
$$ 
where $\varphi_{\omega_C \otimes \eta^i}$ is the map given by the linear system $|\omega_C \otimes \eta^i|$.
\end{lem}

The proof is a slight modification of the proof of \cite[Proposition 12.5.3]{bl}, which we omit. In particular, for $n=2$ 
the composition $\varphi_{\omega_C \otimes \eta} \circ f$ coincides with the projectivized differential of $\pi$.
This is used in \cite[Corollary 12.5.5]{bl} to find conditions for the differential $d\pi_p$ 
to be injective at a point $p \in \tC$ in the case of $s \leq 1$. In any case, we have as a direct consequence of Lemma \ref{lem3.2},

\begin{prop}
The differential of the Abel-Prym map $\pi_c$ is injective at a point $p \in \tilde{C}$, if the point $f(p)$ 
is not a base-point of the linear system $|\omega_C \otimes \eta^i|$ for some $i, \; 1 \leq i \leq n-1$.
\end{prop}

\begin{cor}  \label{cor3.4}
Suppose $n \geq 3$ and either $s \geq 1$ or $s= 0$ and $\Cliff (C) \geq 2$, or $n=2$ and $s \geq 2$.
Then the differential $d\pi_p$ 
of the Abel-Prym map is injective at any point $p \in \tC$.
\end{cor}

\begin{proof}
Suppose first $n \geq 3$. For $s \geq 1$ the line bundle $\omega_C \otimes \eta^2$ is globally generated according to 
Lemma \ref{lem2.1}. For $s=0$, Lemma \ref{lem2.1} gives that $\omega_C \otimes \eta$ is globally generated, unless
$\eta = \cO_C(x-y)$ with $x \neq y \in C$. But in the last case $\omega_C \otimes \eta^2$ is globally generated, since
$\eta^2 = \cO_C(2x - 2y)$ and $2x - 2y \sim u - v$ would imply $\Cliff (C) \leq 1$. For $n=2$ and $s \geq 2$ the line bundle
$\omega_C \otimes \eta$ is globally generated according to Lemma \ref{lem2.1}.
\end{proof}

Combining Proposition \ref{prop3.1} and Corollary \ref{cor3.4} we get the main result of this section.

\begin{thm}
{\em(1)} Suppose $n \geq 3$ and either $s \geq 1$ or $s=0$ and $\Cliff (C) \geq 2$ or $n=2$, $s \geq 2$ 
and $\tC$ not hyperelliptic.
Then the Abel-Prym map $\pi$ is an embedding of $\tC \setminus f^{-1}(B)$ and maps the ramification divisor $f^{-1}(B)$ to an 
ordinary $r$-fold point.\\
{\em (2)} Suppose $n=2, \; s \geq 2$ and $\tC$ hyperelliptic. Then the Abel-Prym map $\pi$ is a double covering 
mapping the ramification divisor $f^{-1}(B)$ to an ordinary $r$-fold point.
\end{thm}

In the missing case $n=2, \; s \leq 1$ and $\tC$ hyperelliptic $\pi: \tC \ra D \subset P$ is a double covering 
onto a smooth curve $D$ such that the Prym variety of $f$ is the Jacobian of $D$ (see \cite[Corollary 12.5.7]{bl}.

\section{The differential of the Prym map}

Let $\mathcal{R}_{g}(n,r)$ denote the coarse moduli space of triples $(C, B, \eta)$, where $C$ is a 
smooth projective curve of genus $g$, $B$ an effective reduced divisor of degree $r \geq 0$ on $C$ and $\eta $ 
a line bundle on $C$ satisfying $\eta^{\otimes n} \simeq\cO_C(B)$. We assume that for $r = 0$, $\eta$ is a 
proper $n$-division point of the Jacobian. Equivalently, $\mathcal{R}_{g}(n,r)$ is the moduli 
space of $n$-fold cyclic coverings $f: \tC \rightarrow C$  with smooth irreducible projective 
curves $C$ of genus $g$ and $\tC$ of genus
$$
g(\tC) = ng -n + 1 + \frac{r(n-1)}{2}.
$$
Note that $f$ is totally ramified over every point of the divisor $B$. 
 The  Prym variety $P = \im (1 - \sigma)$ of the covering $f: \tC \rightarrow C $
is an abelian subvariety of the Jacobian $J\tC$ of dimension  
$$
p=(n-1)(g-1) + \frac{r(n-1)}{2}.
$$
According to \cite[Corollary 12.1.4 and Lemma 12.3.1]{bl} the canonical polarization of $J\tC$ induces a polarization of type
$$
D =  (1, \dots ,1,n, \ldots ,n),
$$
where $1$ occurs $p - (g-1)$ times and $n$ occurs $(g-1)$ times if $r=0$, and  
$1$ occurs $p - g$ times and $n$ occurs $g$ times, if $r > 0$.\\

Let $\mathcal{A}_{p, D}$ 
denote the moduli space of polarised abelian varieties of dimension $p$ with polarisation 
of  type $D$. 
The construction associating to every $n$-fold covering its polarized Prym variety   
defines a map 
$$
Pr = Pr_g(n,r): \mathcal{R}_{g}(n,r) \longrightarrow \mathcal{A}_{p, D},
$$
called the {\it Prym map} (of type $(g,n,r)$). It is a morphism, since the construction works also for  
families of the objects involved. 
In this section we describe the differential of the Prym map at a point 
$(C,B,\eta) = [f: \tC \ra C] \in \mathcal{R}_{g}(n,r)$.
It turns out that it is easier to describe its dual, the codifferential of $Pr$.  \\

Let the notation  be as in Section 3. By Serre duality we have $(H^0(\tC, \omega_{\tilde C}))^* = H^1(\tC, \cO_{\tilde C})$.
Hence the decomposition \eqref{dec} of $H^0(\tC, \omega_{\tilde C})$ induces a decomposition
$$
H^1(\tC, \cO_{\tilde C}) = \bigoplus_{i=0}^{n-1} H^1(C, \eta^{i-n})
$$
into eigenspaces of the action of $G$. Analogously to \eqref{eq3.2} we denote by
$$
H^1(\tC, \cO_{\tilde C})^+ = H^1(C, \cO_{\tilde C}(-B)) \quad \mbox{and} \quad 
H^1(\tC, \cO_{\tilde C})^- = \bigoplus_{i=1}^{n-1} H^1(C, \eta^{i-n})
$$
the eigenspace with eigenvalue 1 and its complement, the sum of the other eigenspaces.

Equation \eqref{eq3.3} implies that the tangent space
$t_P$ of $P$ at the origin is 
\begin{equation} \label{tang}
t_P= (t_{J\tC} )^- = (H^0(\tC, \omega_{\tilde C})^-)^* = H^1(\tC, \mathcal{O}_{\tC})^-,
\end{equation}
which gives for the cotangent space $t_P^*$ of $P$ at the origin,
$$
t_P^*= H^0(\tC, \omega_{\tC})^- = \bigoplus_{i=1}^{n-1} H^0(C, \omega_C \otimes \eta^{n-i}).
$$

Let $\mathcal{M}_{g}(r)$ denote the coarse moduli space of $r$-pointed smooth projective curves of genus $g$.
Since the forgetful map $\mathcal{R}_{g}(n,r) \rightarrow \mathcal{M}_{g}(r)$, \; $(C, B, \eta) \mapsto (C,B)$ is \'etale,  
the tangent space to $\mathcal{R}_{g}(n,r)$ at 
$(C, B,\eta)$ is isomorphic to the tangent space to $\mathcal{M}_{g}(r)$ at $(C,B)$. According to \cite[Example 3.4.19]{s}
or \cite[p. 94]{hm} this space is 
isomorphic to $H^1(C, T_C(-B) )$ and its dual is isomorphic to $H^0(C, \omega_C^{\otimes 2}(B) )$.
 
On the other hand, the cotangent space to $\mathcal{A}_{p,D}$ at the point $P$ can be identified with the second
symmetric product $S^2(H^0(\tC, \omega_{\tC})^-))$. From \eqref{eq3.2} we deduce an isomorphism
$$
S^2 (H^0( \omega_{\tC})^-) \simeq \bigoplus_{i=1}^{n-1} S^2 H^0( \omega_C \otimes \eta^{n-i} )  \oplus  
\bigoplus_{1\leq j < k \leq n-1} H^0(\omega_C\otimes \eta^{\eta -j}) \otimes H^0(\omega_C \otimes \eta^{\eta -k}) .
$$
The eigenspace of $S^2(H^0( \omega_{\tC})^-)$ corresponding to the character $\chi^{\nu}$ of the group 
generated by $\sigma$ is 
\begin{equation*} 
\bigoplus_{\stackrel{1\leq i \leq n-1}{2i \equiv \nu \textrm{ mod } n} } S^2 H^0( \omega_C\otimes \eta^{n-i} )  \oplus  
		\bigoplus_{ \stackrel{1 \leq j < k \leq n}{j+k \equiv \nu \textrm{ mod } n }}  
		H^0(\omega_C\otimes \eta^{n -j}) \otimes H^0( \omega_C \otimes \eta^{n-k}).
\end{equation*}  
We obtain the following commutative diagram compatible with  the action of $ G =\langle \sigma \rangle$ :
\begin{equation}  \label{equivariant}
		\hspace{4cm}  S^2 (H^0( \omega_{\tC})^-)  \hspace{1.2cm}  \longrightarrow \hspace{1.2cm}
		\bigoplus_{\nu=2}^{2n-2} H^0(\omega_C^2 \otimes \eta^{2n-\nu} )
\end{equation} 
$ 
\hspace{6cm} \downarrow \ p \hspace{6.4cm} \downarrow \ p^+ 
$
$$	
		[S^2( \omega_C \otimes \eta^{\frac{n}{2}})]  \oplus	\bigoplus_{ j=1}^{j=[\frac{n}{2}]} H^0(\omega_C   
                   \otimes   
		\eta^{n -j})  \otimes H^0( \omega_C \otimes \eta^{j})  \hspace{.4cm}  \stackrel{\mu} {\longrightarrow}
		\hspace{.4cm}
                   H^0( \omega_C^2(B) )
 $$     
where the factor in square parenthesis  on the bottom row does not occur if $n$ is odd.  The vertical arrows are the 
projection on the $\sigma$-invariant part  and $\mu$ is the multiplication of sections on every factor
of the direct sum. 

\begin{prop} \label{main}
With the identifications above, the codifferential of Prym map at the point $(C,B,\eta) = [f: \tC \ra C] \in \mathcal{R}_{g}(n,r)$ can be identified with the canonical map
$$
\varphi_{C,\eta} :  S^2(H^0(\tC, \omega_{\tC})^-)  \longrightarrow H^0(C, \omega_C^{\otimes 2}(B)),
$$
where $\varphi_{C,\eta}$ is the composed map $\mu \circ p$ of diagram \eqref{equivariant}. 
\end{prop}

In the case $n = 2, \; r = 0$ this was proved \cite[Proposition 7.5]{Be} and stated for arbitrary $n$ and $r = 0$ in 
\cite[Proposition 4.6]{t}. Beauville's proof generalizes directly to this more general case, but seems not to generalize to $r > 0$.  For $n = r = 2$ the proposition was  given in \cite[p. 123]{bcv}.  In order to prove the Proposition we shall need
the following lemmas. For any variety $X$ we denote by $T_X$ its tangent sheaf.

\begin{lem} \label{mult}
Let $X$ be a smooth projective curve  and $A= V / \Lambda$ an abelian variety. Assume there is a non constant morphism  
$\pi : X \rightarrow  A$  whose image generates $A$ as abelian variety. 
Then the dual of the map
$$
H^1(d\pi) : H^1(X, T_X) \lra  H^1(X, \pi^*T_A)
$$
coincides with the multiplication of sections map
$$
V^* \otimes H^0(X, \omega_X) \subset H^0(X, \omega_X) \otimes  H^0(X, \omega_X) \lra H^0(X, \omega^2_X)
$$
with respect to the identifications 
$$
H^1(X, \pi^*T_A)^* \simeq V^* \otimes H^0(X, \omega_X) \quad  \mbox{and} \quad
H^1(X, T_X) ^* \simeq H^0(X, \omega_X^2).
$$
\end{lem}
\begin{proof}
Consider a non constant morphism $\pi : X \ra  A= V / \Lambda$. The differential of $\pi$ gives an  inclusion of sheaves 
$ d\pi: T_X \hra \pi^*T_A = V \otimes \cO_X $ ; its dual map
\begin{equation} \label{map}
 V^* \otimes  \cO_X  \lra  \omega_X 
\end{equation}
has rank 1. Since the image of $X$ generates $A$ as abelian variety, $V^*$ is a subspace of 
$H^0(X, \omega_X )$ and then this map is the  evaluation of sections.  Tensoring  (\ref{map}) by $\omega_X$ and taking  
global sections, we get the multiplication of sections
$$
 V^* \otimes H^0(X, \omega_X) \ra H^0(X, \omega_X^2),
$$
 which is the dual map of $H^1(d\pi)$ .
\end{proof}

According to the definition in Section 3 the Abel-Prym map $\pi: \tC \ra P$ factorizes via the Abel-Jacobi map $\alpha$:
\begin{eqnarray}\label{abel}
\xymatrix{
    		\tC \ar[rr]^{\pi} \ar[dr]_{\alpha}&& P\\
		& J\tC \ar[ur] _{1-\sigma}&
	}	
\end{eqnarray}
The polarization $\Xi$ on $P$ induces an isogeny $\varphi_{\Xi} : P \ra \widehat{P}$ onto the dual abelian variety $\widehat{P}$ 
and thus an isomorphism 
$t_P \ra t_{\widehat{P}}$. According to \eqref{tang} this permits us to identify $t_{\widehat{P}}$ with $ H^1(\tC, \cO_{\tC})^-$. 
\begin{rem}\label{rmk4.3}
This identification is also given by the differential of the dual $\widehat{\pi}$ of the map $\pi$ at the origin:
$$
d\widehat{\pi}: t_{\widehat{P}} \lra t_{J\tC}=H^1(\tC, \cO_{\tC}),
$$
which maps $t_{\widehat{P}}$ isomorphically onto $H^1(\tC, \cO_{\tC})^-$.
\end{rem}

The isomorphism 
$T_{J\tC} \simeq t_{J\tC} \otimes \cO_{J\tC}$ induces an isomorphism 
\begin{equation} \label{eq4.4}
 H^1(J\tC,  T_{J\tC} ) = H^1(\tC, \cO_{\tC})  \otimes  H^1(\tC, \cO_{\tC} ).
\end{equation}
On the other hand, since apart from  $H^1(\tC, \cO_{\tC})^-$, also $H^1(P, \cO_P)$ can be considered as the tangent space of 
$P$ at the origin, we have an equality 
$$
H^1(P,\cO_P) = H^1(\tC, \cO_{\tC})^-.
$$
Hence the isomorphism 
$T_P \simeq t_{\widehat{P}} \otimes \cO_P$ implies that we can identify 
\begin{equation} \label{eq4.5}
 H^1(P,T_P) = H^1(\tC, \cO_{\tC})^- \otimes  H^1(P, \cO_P) = H^1(\tC, \cO_{\tC})^- \otimes H^1(\tC, \cO_{\tC})^-
\end{equation}
as well as
\begin{equation} \label{eq4.6}
 H^1(\tC,  \pi^*T_P )  =  H^1(\tC, \cO_{\tC})^- \otimes  H^1(\tC, \cO_{\tC} ). 
\end{equation}

\begin{rem}\label{rmk4.4}
 The dual of this last identification is provided by the differential of the Prym map $\pi$:
$$
H^0(\tC, \pi^*\Omega_P \otimes \omega_{\tC}) = H^0(\tC, \pi^*\Omega_P) \otimes H^0(\tC, \omega_{\tC}) \stackrel{d\pi \otimes id}{\lra} H^0(\tC, \omega_{\tC})^- \otimes  H^0(\tC, \omega_{\tC} ).
$$
where $d\pi$ induces an isomorphism $ H^0(\tC, \pi^*\Omega_P) \simeq H^0(\tC, \omega_{\tC})^- $.
\end{rem}

\begin{lem} \label{lemma4.3}
 The  diagram 
\begin{eqnarray} \label{diagram}
\xymatrix{ 
		H^1(\tC, T_{\tC} ) \ar[d]_{H^1(d\alpha)} \ar[r]^{H^1(d\pi)} &  H^1(\tC , \pi^* T_P)\\
                H^1(J\tC, T_{J\tC}) \ar[r]^{p^- \otimes p^-} & H^1(P, T_P)\ar[u]^{\pi^*}
         }
\end{eqnarray}
commutes, where $p^-$ is the projection  onto the subspace $H^1(\tC, \cO_{\tC})^-  $  using the identifications above.
\end{lem}

\begin{proof}
Note first that $H^1(d\alpha)$ is a map with target $H^1(\tC, \alpha^*T_{J\tC})$. 
Since $T_{J\tC}$ is a trivial bundle we may identify $H^1(\tC, \alpha^*T_{J\tC}) = H^1(J\tC,T_{J\tC})$ which gives the map
of the left hand vertical arrow of diagram \eqref{diagram}.
By Remark \ref{rmk4.3} , (\ref{eq4.5}) and (\ref{eq4.6})  the map $\pi^*$ corresponds to the map
$$
d\widehat{\pi} \otimes id : H^1(P,\cO_P)\otimes t_P \lra H^1(\tC , \cO_{\tC}) \otimes t_P,
$$
whose image is $H^1(\tC , \cO_{\tC})^- \otimes t_P$.

%---- comment----------------------------------------
%On the other hand, by Serre duality we have the commutative diagram:
%\begin{eqnarray}
%\xymatrix{
%    		 H^1(\tC , \pi^* T_P)^* \ar[r]^{{}^tH^1(d\pi)} \ar[d]_{\simeq} &  H^1(\tC, T_{\tC} )^* \ar[d]^{\simeq} \\
%                 H^0(\tC,\pi^*\Omega_{P}\otimes\omega_{\tC} ) \ar[r]^{ \quad d\pi\otimes id } & H^0(\tC, \omega_{\tC}\otimes\omega_{\tC})
%	}	
%\end{eqnarray}
%-------------------------------------------------

Therefore the diagram \ref{diagram} follows from the commutativity of \ref{abel} and identifications
(\ref{eq4.4}) and (\ref{eq4.5}) .
\end{proof}

\begin{proof} \emph{of Proposition \ref{main}}.
As noted above, the tangent space of $\cR_g(n,r)$ at the point $(C,B,\eta)$ can be identified with $H^1(C,T_C(-B))$.
It coincides with the versal deformation space of the corresponding covering $f: \tC \ra C$. 

%On the other hand, the tangent space $t_P$ of $P$ at $0$ is $H^1(\tC,\cO_{\tC})^-$.

The identification  $t_P = t_{\widehat{P}}$ gives an involution $j$ on $t_P \otimes t_{\widehat{P}}$ 
interchanging the factors. The deformation space of the abelian variety $P$ as a 
complex manifold is $H^1(P,T_P) = t_P \otimes t_{\widehat{P}}$.
The subspace of infinitesimal deformations of $P$ which preserve the polarization $\Xi$ is
$$
(t_{P} \otimes t_{\widehat{P}} )^{j} = S^2(t_P) = S^2(H^1(\tC,\cO_{\tC})^-).
$$
This means that the tangent space of $\cA_{p,D}$ at the point $P$ is $S^2(H^1(\tC,\cO_{\tC})^-)$.

Now the natural cup product map $H^1(\tC, T_{\tC}) \times H^0(\tC, \omega_{\tC}) \ra H^1(\tC, \cO_{\tC})$
is compatible with the action induced by the automorphism $\sigma$ and hence induces a map
$$
H^1(\tC, T_{\tC})^+ \times H^0(\tC, \omega_{\tC})^- \lra H^1(\tC, \cO_{\tC})^-.
$$ 
The image of the corresponding map
\begin{eqnarray*}
H^1(\tC, T_{\tC})^+ \ra &\Hom(H^0(\tC, \omega_{\tC})^-,H^1(\tC, \cO_{\tC})^-) \\ 
& = (H^0(\tC, \omega_{\tC})^-)^* \otimes H^1(\tC, \cO_{\tC})^- \\
& =  H^1(\tC, \cO_{\tC})^- \otimes H^1(\tC, \cO_{\tC})^-
\end{eqnarray*}
is invariant under the involution $j$. Hence it induces a map 
\begin{equation*}
H^1(C,T_C(-B)) = H^1(\tC,T_{\tC})^+ \lra S^2(H^1(\tC,\cO_{\tC})^-), 
\end{equation*}
which is the differential of the Prym map $Pr: \cR_g(n,r) \lra \cA_{p,D}$ at the point $(C,B,\eta)$. 
From Lemma \ref{lemma4.3} we conclude that this map can be considered as a map
$$
H^1(\tC,T_{\tC})^+ \lra H^1(\tC, \pi^*T_P), 
$$
whose image is $S^2 (H^1(\tC,\cO_{\tC})^-) \subset H^1(\tC, \cO_{\tC})^- \otimes  H^1(\tC, \cO_{\tC} ) = H^1(\tC,  \pi^*T_P )$ 
(see \eqref{eq4.6}).
Hence according to Lemma \ref{mult}  the codifferential of $Pr$ at the point $(C,B,\eta)$
is the canonical map 
$$ 
S^2(H^0(\tC,\omega_{\tC})^-) \lra H^0(C,\omega_C^2(B)),
$$
given by multiplication of sections.
\end{proof}

\section{Injectivity of the differential of the Prym map}

Proposition \ref{main} implies that in order to show that the differential of the Prym map 
$Pr : \mathcal{R}_{g}(n,r) \lra \mathcal{A}_{p, D}$ at a point $(C,B,\eta)$ is injective, it suffices to show that there is an
$i,\; 1 \leq i \leq n-1$ such that the canonical map
$$
\mu: H^0(\omega_C \otimes \eta^i) \otimes H^0(\omega_C \otimes \eta^{n-i}) \lra H^0(\omega_C^2(B))
$$
is surjective.

First we consider the case of an \'etale covering, i.e. $r=0$. We denote the moduli space of $n$-fold cyclic 
\'etale coverings of curves of genus $g$ by $\mathcal{R}_{g}[n]$.
In the case of $n$ even we deduce:

\begin{prop} \label{prop5.1}
Suppose $\Cliff(C) \geq 3$. Then for $n$ even the differential of Prym map 
$Pr : \mathcal{R}_g[n] \longrightarrow \mathcal{A}_{p, D}$ is injective at the point $[f: \tilde{C} \ra C]$.
\end{prop}
 
\begin{proof}
For $n=2m$ one of the  factors of  $S^2(H^0(\tC, \omega_{\tC})^-) $ is 
$S^2H^0(C,\omega_C\otimes \eta^m)$. 
Since $\eta^m$ is a 2-torsion point, the assertion is a consequence of Corollary \ref{cor2.4} and Proposition \ref{main}.
\end{proof}

For arbitrary $n$ we only have:
\begin{prop} \label{prop5.2}
Suppose $\Cliff(C) \geq 2n-1$. Then the differential of Prym map 
$Pr : \mathcal{R}_g[n] \lra \mathcal{A}_{p, D}$ is injective at the point $[f: \tilde{C} \ra C]$.
\end{prop}
\begin{proof}
Again this is a consequence of the remarks before Proposition \ref{prop5.1}, Corollary \ref{cor2.4} 
and Proposition \ref{main}.
\end{proof}

\begin{rem}
As we saw in the proof of Corollary \ref{cor2.4}, the assumptions on the Clifford index imply $g \geq 8$ in 
Proposition \ref{prop5.1} and $g \geq 4n-1$ in Proposition \ref{prop5.2}. 
Hence the bound for the Clifford index is not sharp in Proposition \ref{prop5.1}, since it is also valid for $g=6, n=2$ 
(see \cite[Corollaire 7.11]{Be}).

The bounds in Proposition \ref{prop5.2} are certainly not sharp if $n$ is a composed number. If $p$ is a proper prime divisor of $n$, then $\eta^{\frac{n}{p}}$
is a $p$-division point and the method of Proposition \ref{prop5.1} gives a better result. We omit this, 
since it is easy to work out.
\end{rem}

For a general covering in $\cR_g[n]$ there is a better bound which relies on the following lemma.

\begin{lem} \label{lem5.4}
For $g\geq 5$ there exists a curve $C$ of genus $g$ admitting 
a proper $n$-division point $\eta \in JC [n]$ such that the line bundle
$\omega_C \otimes \eta$ is very ample for any $n \geq 2$.
\end{lem}
\begin{proof}
Consider the locus 
$$ 
\mathcal{V}_n := \{([C], \eta)  \in \mathcal{R}_g[n] \mid  \ \omega_C \otimes \eta  \textrm{ is not very ample}\}.
$$
which is a Zariski-closed subset of $\mathcal{R}_g[n]$. According to Lemma \ref{lem2.1}, a curve $C$ such that $([C], \eta) 
\in \mathcal{V} $ admits points $x,y,u,v$ with $x+y \neq u + v$ such that $nx + ny \sim nu + nv$. 
According to the Riemann-Hurwitz formula the corresponding $2n:1$ covering $C \ra \PP^1$ 
is branched 
at most over 
$$
b = 4n + 2g -2 -4(n-1) +2 = 2g + 4
$$
points of $\PP^1$.
Hence, the dimension  of $ \mathcal{V}_n$ is at most the dimension of Hurwitz scheme
$$
\mathcal{H}_{2n,b} =  \{ \pi: C \stackrel{2n:1}{\lra} \mathbb{P}^1 \textrm{ branched over $b$ points } \}.
$$
It is known that the Hurwitz schemes have at most the 
expected dimension $b$, the number of branch points of the cover.  Therefore, the dimension 
of locus in $\mathcal{M}_g$ of curves $C$ admitting a $g^1_{2n}$ of the form $| n x + n y|$ is at 
most  
$$
\dim \mathcal{H}_{2n,b} - \dim \Aut (\mathbb{P}^1) = 2g+1, 
$$ 
which is less than $\dim \mathcal{R}_g[n] = \dim \mathcal{M}_g =3g-3$ for $g \geq 5$.  This completes the proof. 
\end{proof}

Using this, we can conclude,

 \begin{prop} \label{prop5.5}
For any $n\geq 2$ and  $g \geq 7$  
the  Prym map $Pr :  \mathcal{R}_g[n] \lra \mathcal{A}_{p,D}$ is generically finite.
In particular the image of $Pr$ is of dimension $3g-3$. 
\end{prop}

\begin{proof}
It suffices to show that there exists an element $([C], \eta) \in \cR_g[n]$ such that the differential of the Prym map is 
injective at $([C], \eta)$. For this it suffices to show, according to Proposition \ref{main}, that
there exists a $([C], \eta) \in \cR_g[n]$ such that the canonical map 
$H^0(C, \omega_C \otimes \eta) \otimes H^0(C, \omega_C \otimes \eta^{-1}) \ra H^0(C, \omega^2)$
is surjective.

This follows from \cite[Theorem 1]{b} using Lemma \ref{lem5.4} and the fact that a general curve of genus $g \geq 7$ 
is of Clifford index $\geq 3$.
\end{proof}

Finally, we consider the case of ramified $n$-fold cyclic coverings, i.e. $r > 0$. Here we have the 
following consequences of Proposition \ref{main}:

\begin{prop} \label{prop5.6}
Suppose  $g\geq 2$ and $r \geq 6$ if $n$ is even, repectively $r \geq 7$ if $n$ is odd.
Then the differential of Prym map $Pr : \mathcal{R}_g(n,r) 
\lra \mathcal{A}_{p, D}$ is injective at any point $(C,B,\eta) \in \mathcal{R}_g(n,r)$.
\end{prop}

\begin{proof}
As we have said at the beginning of this section, it suffices to show that the multiplication map
$$
H^0(C, \omega_C \otimes \eta^{[\frac{n}{2}]}) \otimes H^0(C, \omega_C \otimes \eta^{n-[\frac{n}{2}]}) \longrightarrow
H^0(C, \omega_C^2 (B))
$$
is surjective. According to \cite[Theorem 1]{b} this map is surjective if the line bundles 
$\omega_C \otimes \eta^{[\frac{n}{2}]}$ and $\omega_C \otimes \eta^{n-[\frac{n}{2}]}$ are very ample .

Since $\deg (\omega_C ) =2g-2 $, the very ampleness condition is satisfied if 
$\deg \eta^{[\frac{n}{2}]} $ (resp. $\deg \eta^{n-[\frac{n}{2}]} $) is at least 3.  
For $n=2m$  even, $[\frac{n}{2}] -  [\frac{n}{2}] =m $ and then 
$\deg \eta^m= \frac{mr}{n} = \frac{r}{2} \geq 3 $ if and only if  $r \geq 6$.
On the other hand, if $n= 2m+1$, $m\geq 1$, the line bundles 
$\omega_C \otimes \eta^{m}$ and $ \omega_C \otimes \eta^{m+1}$ 
are very ample as soon as 
$ \frac{r}{n} (\frac{n-1}{2}) \geq 3$.
This is the case for $r \geq 7$, since $r$ is a multiple of $n$ according to Hurwitz formula.
\end{proof}

Since $n$ divides $r$, we are left with 6 cases for the pair $(n,r)$, for which we have a similar result. It is given in the following corollary, 
the proof again is an application of Lemma \ref{lem2.1} and \cite[Theorem 1]{b}. We omit the details. 

\begin{prop} \label{prop5.7}
The differential of Prym map $Pr : \mathcal{R}_g(n,r) 
\lra \mathcal{A}_{p, D}$ is injective at the point $(C,B,\eta) \in \mathcal{R}_g(n,r)$ in the following cases:
\begin{itemize}
\item for $n=2, r=2$ i.e. $\deg \eta =1$, if $\eta \neq \cO_C(x+y-u)$ and $\Cliff(C) \geq 2$;
\item for $n=2, r=4$ i.e. $\deg \eta =2$, if $\eta \neq \cO_C(x+y)$ and $\Cliff(C) \geq 1$;
\item for $n=3, r=3$ i.e. $\deg \eta =1$, if $\eta \neq \cO_C(x+y-u)$ and $\Cliff(C) \geq 2$;
\item for $n=3, r=6$ i.e. $\deg \eta =2$, if $\eta \neq \cO_C(x+y)$ and $\Cliff(C) \geq 1$;
\item for $n=4, r=4$ i.e. $\deg \eta =1$, if $\eta^2 \neq \cO_C(x+y)$ and $\Cliff(C) \geq 1$;
\item for $n=5, r=5$ i.e. $\deg \eta =1$, if $\eta^2 \neq \cO_C(x+y)$ and $\Cliff(C) \geq 1$.
\end{itemize}
\end{prop}
 
For $g \geq 3$ respectively $g \geq 5$ there are curves of Clifford index $\geq 1$ respectively $\geq 2$ we get 
as an immediate consequence of Propositions \ref{prop5.6} and \ref{prop5.7},
 
\begin{cor} \label{cor5.8}
Suppose one of the following conditions is satisfied
\begin{itemize}
\item $g \geq 2$ and $r \geq 6$ for $n$ even or $r\geq 7$ for $n$ odd, 
\item $g \geq 3$ and $n=r = 4$ or $5$ or $(n,r) = (2,4)\; or\; (3,6)$,
\item $g \geq 5$ and $n=r = 2$ or 3.
\end{itemize} 
Then  the  Prym map $Pr :  \mathcal{R}_g(n,r) \lra \mathcal{A}_{p,D}  $    is generically finite.
In particular the image of $Pr$ is of dimension $3g-3 + r$. 
\end{cor}


\begin{thebibliography}{999999}

\bibitem[BCV]{bcv}
    F. Bardelli, C. Ciliberto, A. Verra:
    \textit{Curves of minimal genus on a general abelian variety}.
    Compos. Mathem. 96 (1995), 115-147.

\bibitem[Be]{Be} 
    A. Beauville: \textit{Vari\'et\'es de Prym et Jacobiennes intermediares}. 
    Annales Ec. Norm. Sup. 3 (1977),  309-391.

\bibitem[BL]{bl}
    Ch. Birkenhake, H. Lange:
    \textit{Complex Abelian Varieties}. Second edition,
    Grundlehren der Math. Wiss. 302, Springer - Verlag (2004).

\bibitem[Bu]{b}
    D. C. Butler: 
    \textit{Global sections and tensor products of line bundles over a curve}
    Math. Z. 231 (1999), 397-407.
    
\bibitem[D]{d}
    R. Donagi:
    \textit{The tetragonal construction}.
    Bull. Am. Soc. 4 (1981), 181-185.
    
\bibitem[FS]{fs}
    R. Friedman, R. Smith:
    \textit{The generic Torelli theorem for the Prym map}.
    Invent. Math. 67 (1982), 473-490.   
    
\bibitem[GL]{gl} 
     M. Green, R. Lazarsfeld: 
     \textit{On the projectivity normality of complete linear series on an algebraic curve}. 
     Invent.  Math. 83 (1986), 73-90. 
     
\bibitem[HM]{hm}
     J. Harris, I. Morrison:
     \textit{Moduli of curves}.
     GTM, No. 187, Springer - Verlag (1998).
     
\bibitem[K]{k}
     V. Kanev:
     \textit{The global Torelli theorem for Prym varieties at a generic point}.
     Math. USSR-Izv. 20 (1983), 235-258.  
     
\bibitem[LS]{ls} 
     H. Lange, E. Sernesi:
     \textit{On the Hilbert scheme of a Prym variety}.
     Ann. di Matem. 183 (2004), 375-386.       
     
\bibitem[S]{s}
     E. Sernesi:
     \textit{Deformations of algebraic schemes}. 
     Grundlehren der Math. Wiss. 302, Springer - Verlag (2006).
     
\bibitem[T]{t}
     A. Tamagawa:
     \textit{Finiteness of isomorphism classes of curves in positive characteristic with prescribed 
     fundamental groups}. J. Alg. Geom. 13 (2004), 675-724. 
    
\end{thebibliography}
\end{document}